\documentclass{article}
\usepackage[title]{appendix}
\usepackage{amsmath,amssymb,amsthm}
\usepackage{authblk}
\usepackage{color}
\usepackage{graphicx}
\usepackage{float}
\usepackage{pdfpages}
\usepackage{hyperref}
\usepackage{fancyhdr}
\hypersetup{colorlinks=true,linkcolor=blue,citecolor=blue}
\usepackage{lipsum}
\usepackage{cleveref}
\makeatletter
\renewcommand*\l@section{\@dottedtocline{1}{1.5em}{2.3em}}
\makeatother

\crefname{equation}{equation}{equations}
\Crefname{equation}{Equation}{Equations}
\crefrangelabelformat{equation}{(#3#1#4-#5#2#6)}

\crefmultiformat{equation}{equations (#2#1#3}{, #2#1#3)}{#2#1#3}{#2#1#3}
\Crefmultiformat{equation}{Equations (#2#1#3}{, #2#1#3)}{#2#1#3}{#2#1#3}
\newtheorem{theorem}{Theorem}

\newtheorem{proposition}[theorem]{Proposition}
\newtheorem{lemma}[theorem]{Lemma}

\newtheorem{remark}{Remark}
\newtheorem{maintheorem}{Main Theorem}

\author[1]{Daniel Peralta-Salas\thanks{dperalta@icmat.es}}
\author[2]{Miguel Vaquero\thanks{mvaquero@faculty.ie.edu}}
\affil[1]{Instituto de Ciencias Matem\'aticas, Consejo Superior de Investigaciones Cient\'ificas\\
 28049 Madrid, Spain}
\affil[2]{School of Science and Technology, IE University\\ 40003 Segovia, Spain}

\title{Beltrami Fields with Morse Proportionality~Factor}
\date{}
\begin{document}
\maketitle

\begin{abstract}
In this work we study Beltrami fields with non-constant proportionality factor on $\mathbb{R}^3$. More precisely, we analyze the existence of vector fields $X$  satisfying the equations 
$ curl(X)=fX$ and  $div(X)=0$ for a given $f\in C^\infty(\mathbb R^3)$ in a neighborhood of a point $p\in\mathbb{R}^3$. Since the regular case has been treated previously, we focus on the case where $p$ is a non-degenerate critical point of $f$. We prove that for a generic Morse function $f$, the only solution is the trivial one $X\equiv 0$ (here generic refers to explicit arithmetic properties of the eigenvalues of the Hessian of $f$ at $p$). Our results stem from the introduction of algebraic obstructions, which are discussed in detail throughout the paper. 
\end{abstract}


\newpage

\section{Introduction}
A {\it Beltrami field}  in
$\mathbb{R}^3$ is a vector field $X$ satisfying
\begin{align}\label{beltrami}
curl(X)=fX,\qquad div(X)=0,
\end{align}
for some smooth function $f$. When $f$ is a constant, the solutions to Equations~\eqref{beltrami} are sometimes called {\it strong Beltrami fields}. It is well-known that Beltrami fields are 
stationary solutions of the incompressible Euler equations in $\mathbb R^3$. Moreover, they have proven to be 
very powerful tools to analyze the structure
of the solutions to the time-dependent Euler and Navier-Stokes equations, see e.g. \cite{DeLSz,EnLuPe,EnPeLines,EnPeTubes}. In the context of plasma physics, Beltrami fields are known as force--free fields, and define a particularly remarkable class of magnetohydrostatic equilibria. Nonetheless, the
analysis of Beltrami fields with non-constant proportionality factor is extremely
hard. In that regard, one of the main questions, sometimes called the {\it helical flow problem}~\cite{MoYuZa},  is
to determine for which functions $f$ there is a non-trivial vector
field satisfying the Equations~\eqref{beltrami}.

A major result in this direction has been  obtained by Enciso and Peralta-Salas in~\cite{EnPeARMA}, where the authors introduce an operator
$P[f]$ whose vanishing is a necessary condition for the existence of non-trivial solutions to the Equations~\eqref{beltrami}. The
aforementioned operator is constructed in coordinates
adapted to the level sets of the function $f$, through the use of the implicit function
theorem. Remarkably, the
operator $P[f]$ is constructed in an open set $U\subset \mathbb{R}^3$,
maybe smaller than the initial analyzed set, which does not include
critical points. This observation allowed the authors to show that, generically, 
there are no non-trivial solutions to problem~\eqref{beltrami} around regular points, where the gradient of $f$ does not vanish.  The same operator was also used in~\cite{Ab} to investigate the rigidity of Beltrami vector fields.

Following a similar research direction, the authors in~\cite{CeKloARMA} employ Cartan moving frames and exterior differential systems techniques to strengthen the results in~\cite{EnPeARMA}. Notably, all these references~\cite{Ab,CeKloARMA,EnPeARMA} rely on the non-vanishing of the gradient of the function $f$ in the open set under consideration. This restriction leads to the following natural question:
 \begin{center} 
 {\it Do there exist non-trivial solutions to Equations~\eqref{beltrami} in a neighborhood of a critical point of the factor $f$? }
\end{center}


In Appendix~\ref{nonisolated} we provide an example of a function $f$ that has a curve of critical points and admits a non-trivial Beltrami field. However, in this article we shall focus on the generic case where the critical points of $f$ are non-degenerate. Roughly speaking, our first main result shows that in a neighborhood of a non-degenerate critical point of $f$,
the only solution to Equations~\eqref{beltrami} is the trivial one, except for some special cases. More precisely, we prove:
\begin{maintheorem} Assume that $f$ has a non-degenerate critical point at $p$, and the Hessian $d^2f(p)$ satisfies one of the following assumptions:
  \begin{itemize}
  \item The spectrum of $d^2f(p)$ is different from  $\{\alpha,-\alpha,\beta\}$ where $\alpha, \beta\in \mathbb{R}\char`\\ \{0\}$.
    \item The sum of the eigenvalues  of $d^2f(p)$ (the trace) is non-zero.
\item   The spectrum of $d^2f(p)$ is different from  $\{\alpha,\alpha,-i\alpha\}$ with $\alpha\in \mathbb{R}\char`\\ \{0\}$ and $i$ is a natural number with $i\geq 3$.
\end{itemize}
  Then, the only solution to~\eqref{beltrami} in a neighborhood of $p$ is $X\equiv 0$. Obviously, each one of the assumptions is generic.
\end{maintheorem}

The result mentioned above involves the Hessian matrix as an obstruction for the existence of non-trivial Beltrami fields with non-constant factor. Our second main contribution permits the Hessian to be any non-singular matrix, provided that the third and fifth-order terms of the Taylor series expansion vanish.

\begin{maintheorem} Let $f$ be a function with a non-degenerate critical point at $p$.  Assume that the Taylor series expansion of $f$ at $p$ takes the form $f=f_0+f_2+f_4+\mathcal{O}(6)$. That is, the homogeneous terms of order $3$ and $5$ of the function $f$ vanish. Then the unique solution to Equations~\eqref{beltrami} in a neighborhood of $p$ is the trivial one.
\end{maintheorem}

This article leaves open the question of the existence of non-trivial solutions in a neighborhood of a non-degenerate critical point for which none of the main theorems above can be applied. This case is briefly discussed in Appendix~\ref{app:numerical}.

The paper is organized as follows. In Section~\ref{S.stra} we explain the general strategy to prove the two main theorems. The proofs of these theorems are presented in Section~\ref{main}, up to some technical lemmas. To this end, computationally intensive proofs of intermediate results are relegated to the appendices. Finally, in Section~\ref{conclusions} we discuss the case of perturbations of the setting considered in the second main theorem. In the appendices we provide proofs of auxiliary propositions that are instrumental in the proofs of the main theorems. We also provide a brief discussion of the cases that are not covered by the second main theorem, motivated by an example presented in Appendix~\ref{app:numerical}. 

\subsection{Notation}
In this paper we work in the smooth category of $\mathcal{C}^\infty$ functions. We will use $\langle \cdot{}, \cdot{}\rangle$ to denote the scalar product in $\mathbb{R}^3$. Given a function $f$, its differential at point $p$ is denoted by $df(p)$ and its Hessian by $d^2f(p)$. The imaginary unit is denoted by $I$, where $I^2 = -1$. The set of all non-zero real numbers is denoted by $\mathbb{R}\backslash\{0\}$. The real and imaginary part of a complex number, $z$, is $Re(z)$ and $Im(z)$.


\section{Proposed Strategy}\label{S.stra}
Let $f:\mathbb{R}^3\rightarrow \mathbb{R}$ be a smooth function, i.e., in $\mathcal{C}^{\infty}(\mathbb{R}^3)$. We assume the existence of a point $p\in\mathbb{R}^3$ which is a non-degenerate critical point of $f$. This means the differential vanishes, $df(p) = 0$, and the Hessian $d^2f(p)$ is a non-singular matrix at the point under consideration.

It follows from elliptic regularity~\cite{EnPeARMA} that any solution to Equations~\eqref{beltrami} with $f\in \mathcal{C}^{\infty}$ is also $\mathcal{C}^{\infty}$.
To investigate whether the Equations~\eqref{beltrami} have non-trivial solutions in a neighborhood of $p$, we rely on a Taylor series expansion of both the vector field $X$ and the function $f$ in Equations~\eqref{beltrami} around the point $p$. We write $X=X_0+X_1+X_2+\ldots$ where $X_i$ denotes the homogeneous component of degree $i$, and  proceed similarly with the proportionality factor
$f=f_0+f_2+f_3+\ldots$ Notice that $f_1=0$ as $p$ is  a
critical point. We substitute the above series expansion into the Equations~\eqref{beltrami}, which yields
\begin{subequations}\label{expansion}
\begin{align}
&curl(X_0+X_1+X_2+\ldots)=(f_0+f_2+f_3+\ldots)(X_0+X_1+X_2+\ldots),  \label{expansion:1}\\
&  div(X_0+X_1+X_2+\ldots)=0, \label{expansion:2}\\
&  \langle \nabla (f_0 + f_2 + f_3 + \ldots),X_0 + X_1 + X_2 + \ldots \rangle = 0.  \label{expansion:3}
\end{align}
\end{subequations}
Observe that we added the redundant equation $\langle\nabla f, X\rangle = 0$, resulting in Equation~\eqref{expansion:3}. This can be derived from Equations~\eqref{beltrami} by taking the divergence on both sides of the first equation in~\eqref{beltrami} and noting that $X$ is divergence-free. It is important to notice that the least order monomial in~\eqref{expansion:3} is $\langle \nabla f_2,X_0 \rangle$ and, as a consequence of the non-degeneracy of the Hessian of $f$ at $p$, we obtain 
$$X(p) = X_0 = 0\,.$$

In order to investigate the solutions of Equations~\eqref{expansion}, monomials of the same degree can be matched, resulting in an infinite-dimensional
linear system  of equations. Nonetheless, the complexity of this infinite-dimensional system hinders the attainment of results. Our main contributions are based on the observation that  finite-dimensional systems can be extracted from Equations~\eqref{expansion}, which provide enough information to demonstrate that the only solution to the Equations~\eqref{expansion} in a neighborhood of $p$ is the trivial one.

As shown in Proposition~\ref{con2}, if there were  solutions $X\not\equiv 0$ to~\eqref{expansion}, then there should be a first
non-trivial term in this series $\sum_{i=0}^\infty X_i$, say $X_{i_0}$. This term satisfies, among other constraints, the following system of equations
\begin{equation}\label{nose1}
\begin{aligned}
&curl(X_{i_0})=0,\\
&div(X_{i_0})=0,\\
&\langle \nabla f_2, X_{i_0} \rangle=0.
\end{aligned}
\end{equation}
After a careful analysis, we observe that  Equations~\eqref{expansion} only have a non-trivial solution $X_{i_0}$ if a strong constraint on the eigenvalues of the Hessian matrix $d^2f(p)$ is satisfied, as described in Proposition~\ref{con2}. Thus, Proposition~\ref{con2} determines large families of functions $f$ for which
all the terms $X_i$ have to vanish in the series
expansion of $X$. Accordingly, for the aforementioned families of functions, the only possible solution in a neighborhood of $p$ is the trivial one $X \equiv 0$ (see Appendix~\ref{app:strong_continuation}). 

In the case where the function $f$ is such that Equations~\eqref{nose1}
admit a non-trivial solution for some index $i_0$, additional
  equations taken from~\eqref{expansion} can be analyzed. Under suitable generic conditions, this analysis leads to the conclusion that the only compatible solution is $X_i \equiv 0$.
 Following this reasoning, we are able to prove that $X\equiv 0$ is the unique solution for all the functions $f$ that fall
 into one of the categories studied in this paper.

\section{Proof of the Main Theorems}\label{main}

Let $f:\mathbb{R}^3\rightarrow \mathbb{R}$ be a smooth function and
$p\in\mathbb{R}^3$ a non-degenerate critical point of $f$; we can safely assume that $p=0$. We assume that $d^2f(p)$ is a  diagonal matrix of the form
\[
\left(\begin{array}{ccc}2 \sigma_1 &0&0\\0&2\sigma_2&0\\ 0&0&2\sigma_3\end{array}\right),
\]
where $\sigma_1,\ \sigma_2$ and $\sigma_3$ are different from zero, by the non-degeneracy
condition. If that is not the case, we can apply an isometry to transform the Hessian into the desired diagonal form. It is important to note that isometries commute with the $curl$ and the $div$ operators. The inclusion of the superfluous $2$ is solely to ensure that the term $f_2$ equals $\sigma_1x^2 + \sigma_2y^2 + \sigma_3z^2$.

 We decompose the vector field $X$
and the function $f$ into their homogeneous components (Taylor series) around the point $p$
under consideration
$X=X_0+X_1+X_2\ldots$ and $f = f_0 + f_2 + f_3 \ldots$ 
When the obtained expressions are replaced into Equations~\eqref{beltrami}, if the vector field $X$ is not identically zero, then there is a first non-trivial homogeneous polynomial vector field $X_{i_0}$ satisfying the Equations~\eqref{nose1}.  One of the main observations of this paper (Proposition~\ref{con2} below), shows that these algebraic equations only
have a non-trivial solution $X_{i_0}$, when the degree of the polynomial
vector field $i_0$ and
the eigenvalues of $d^2f(p)$ are related in a very specific way. Before proving the main theorems, we introduce some preliminary results.

\begin{lemma}\label{lemma2} If  $\sigma_1, \ \sigma_2$
  and $\sigma_3$ all have the same sign,
  then the only homogeneous polynomial solution to Equations~\eqref{nose1} is $X_{i_0}\equiv0$.
\end{lemma}

\noindent\textbf{Proof:}
Let us prove a slightly stronger result. Let $Y$ be a vector field, not necessarily polynomial. Then,  if the real numbers $\sigma_1, \ \sigma_2$
  and $\sigma_3$ have the same sign, and $Y$ satisfies
\begin{align*}
&curl(Y)=0, \\
&\langle \nabla(\sigma_1 x^2+\sigma_2 y^2+\sigma_3 z^2), Y\rangle =0,
\end{align*}
we claim that $Y\equiv 0$. Indeed, since $curl(Y)=0$, 
then $Y=\nabla g$ for some function $g$ on $\mathbb{R}^3$. Using $\langle
\nabla(\sigma_1 x^2+\sigma_2 y^2+\sigma_3 z^2),Y\rangle = \langle
\nabla(\sigma_1 x^2+\sigma_2 y^2+\sigma_3 z^2),\nabla g\rangle  =0$, we infer that $g$ is a first integral of
the linear vector field $\nabla(\sigma_1 x^2+\sigma_2 y^2+\sigma_3 z^2) $. But since
$\nabla(\sigma_1 x^2+\sigma_2 y^2+\sigma_3 z^2)$ is a linear vector field where all the eigenvalues have
the same sign, the origin is then either a source or a sink. Given that the only continuous
first integrals are constant functions, it follows that  $g=constant$, and consequently, $Y =\nabla g\equiv 0$, as claimed.
\qed
\vspace*{0.25cm}

As a straightforward consequence we obtain the following proposition.

\begin{proposition}\label{samesign} If the eigenvalues of $d^2f(p)$ have the same sign, then the unique solution to Equations~\eqref{beltrami} is $X\equiv 0$.
\end{proposition}
\noindent\textbf{Proof:}  We proceed by induction. First, as previously observed, $X_0=0$. Assume that
$X_{i}\equiv 0$ for $i=0,1,\ldots,n$ and let us show that $X_{n+1}\equiv 0$. By
the induction hypothesis, it is easy to see that $X_{n+1}$ has to satisfy 
\begin{align*}
&curl(X_{n+1})=0,\\
&\langle \nabla(\sigma_1 x^2+\sigma_2 y^2+\sigma_3 z^2), X_{n+1}\rangle =0.
\end{align*}
An application of Lemma~\ref{lemma2} combined with the strong unique continuation property (see Appendix~\ref{app:strong_continuation}) gives the result.
\qed

\begin{remark} Proposition~\ref{samesign} is related to, but weaker than, Theorem $1.2$ in~\cite{EnPeARMA}.
\end{remark}

The next proposition establishes a relation between the
eigenvalues of $d^2f(p)$  and the $X_i$'s. Specifically, the
eigenvalues of $d^2f(p)$ are going to determine the first possible non-trivial term $X_{i_0}$ in the Taylor series expansion $X=X_0+X_1+X_2+\ldots$ Recall that the spectrum of the Hessian $d^2f(p)$ is given by
$\{2\sigma_1,2\sigma_2,2\sigma_3\}$, and a necessary condition for the existence of non-trivial solutions to Equations~\eqref{beltrami} is that not all eigenvalues have the same sign (see Proposition~\ref{samesign}).

\begin{proposition}\label{con2} Let $X_i$ be a vector field whose components are
  homogeneous polynomials of degree $i$.  Consider the system of equations
 \begin{equation}\label{equation_new}
\begin{aligned}
&curl(X_{i})=0,\\
&div(X_{i})=0,\\ 
&\langle \nabla(\sigma_1 x^2+\sigma_2 y^2+\sigma_3 z^2), X_{i} \rangle=0,
\end{aligned}
\end{equation}
where the unknowns are the vector field $X_i$ and the
non-zero real numbers  $\sigma_1,\ \sigma_2, \ \sigma_3$ which do not all have the same sign. Then, a necessary condition for the existence of solutions $X_i\not\equiv 0$ is that
\begin{itemize}
\item If $i=1$ then $\{\sigma_1,\sigma_2,\sigma_3\}=\{\alpha,-\alpha,\beta\}$ for some non-zero real numbers $\alpha, \beta$.
\item If $i=2$ then $\sigma_1+\sigma_2+\sigma_3=0$.
\item If $i\geq 3$ then $\{\sigma_1,\ \sigma_2, \ \sigma_3\}=\{\alpha,\alpha,-i\alpha\}$ for some non-zero real number $\alpha$. 
\end{itemize} 
Moreover, for
 $i\geq 3$ if  we assume that $\sigma_1=\sigma_2=\alpha$ and
 $\sigma_3=-i\alpha$ the solution is  explicitly given
 by \[X_i=\nabla(p(x,y)\cdot{}z)\,,\] where $p(x,y)$ is a
 homogeneous harmonic polynomial of degree $i$ in $\mathbb{R}^2$. 
\end{proposition}

\noindent\textbf{Proof:} See Appendix~\ref{app:con2}.

\begin{remark}
Notice that if there are no resonances among the eigenvalues, i.e., there are no solutions to  $\sigma_1 k_1+\sigma_2 k_2+\sigma_3 k_3=0$ for $k_1,\ k_2$ y $k_3 \in\mathbb{Z}$, then the
  unique solution to Equations~\eqref{beltrami} is $X\equiv0$.
\end{remark}

The next result is the first main contribution of this article. It mainly states that when the spectrum of $d^2f(p)$ does not fall into one of the categories described in Proposition~\ref{con2}, the only solution to Equations~\eqref{beltrami} is $X\equiv0$. This assertion is sufficient to conclude that for a generic Morse factor  $f$ (in the sense of an open and dense set in the $C^\infty$ topology), the only solution to~\eqref{beltrami} is the trivial one.

\begin{theorem}[Main Theorem 1]\label{th:2} Assume that $f$ has a non-degenerate critical point at $p$, and the Hessian matrix $df^2(p)$ does not fall into one of the following categories:
  \begin{itemize}
  \item The spectrum of $d^2f(p)$ is of the form $\{\alpha,-\alpha,\beta\}$ for some $\alpha, \beta,\in \mathbb{R}\char`\\ \{0\}$.
    \item The sum of the eigenvalues  of $d^2f(p)$ (the trace) is equal to zero.
\item   The spectrum of $d^2f(p)$ is of the form  $\{\alpha,\alpha,-i\alpha\}$ for some $\alpha\in \mathbb{R}\char`\\ \{0\}$ and $i$ a natural number with $i\geq 3$.
\end{itemize}
  Then, the only solution to Equations~\eqref{beltrami} is $X\equiv0$. 
\end{theorem}
\noindent\textbf{Proof:} We proceed by induction. Remember that $X_0 = 0$, and assume that
$X_i \equiv 0$ for $i = 0, 1, \ldots, n$. By the induction
hypothesis, it is evident that $X_{n+1}$ must satisfy Equations~\eqref{equation_new}. It then follows from Proposition~\ref{con2} that $X_{n+1} \equiv 0$, confirming that $X$ must vanish, as we wanted to show.

\qed


The remainder of this section addresses the case where the Hessian matrix of the factor $f$ at $p$ falls into one of the scenarios described in Theorem~\ref{th:2}. For $i\geq 3$, Proposition~\ref{con2} implies that, subject to reparametrization and change of variables, the only solution to Equations~\eqref{equation_new} is $X_i=\nabla(p(x,y)\cdot{}z)$, where $p(x,y)$ is a homogeneous harmonic polynomial of degree $i$. It is well-known that $Re((x+Iy)^i)$ and  $Im((x+Iy)^i)$ form a basis for the homogeneous harmonic polynomials of degree $i$ in the plane. Then,
\[
p(x,y)=\lambda_1 Re((x+Iy)^i)+\lambda_2 Im((x+Iy)^i),
\]
with $\lambda_1,\lambda_2$ real constants. A straightforward computation yields
\begin{align*}
  Re((x+Iy)^i)=\sum_{k=0}^i\left(\begin{array}{c}i\\ k\end{array}\right)\cos((i-k)\pi/2)x^{k}y^{i-k},
  \\
  Im((x+Iy)^i)=\sum_{k=0}^i\left(\begin{array}{c}i\\ k\end{array}\right)\sin((i-k)\pi/2)x^{k}y^{i-k}.
\end{align*}
Finally, we introduce the notation 
\[
X_i^1:=\nabla(Re((x+Iy)^i)\cdot{}z)=\left(\begin{array}{c}\sum_{k=1}^i\left(\begin{array}{c}i\\ k\end{array}\right)k\cos((i-k)\pi/2)x^{k-1}y^{i-k}z \\\noalign{\medskip}
\sum_{k=0}^{i-1}\left(\begin{array}{c}i\\ k\end{array}\right)(i-k)\cos((i-k)\pi/2)x^{k}y^{i-k-1}z\\\noalign{\medskip}
\sum_{k=0}^i\left(\begin{array}{c}i\\ k\end{array}\right)\cos((i-k)\pi/2)x^{k}y^{i-k}
\end{array}\right),
\]
and
\[
  X_i^2:=\nabla(Im((x+Iy)^i)\cdot{}z)=
\left(\begin{array}{c}\sum_{k=1}^i\left(\begin{array}{c}i\\ k\end{array}\right)k\sin((i-k)\pi/2)x^{k-1}y^{i-k}z \\\noalign{\medskip}
\sum_{k=0}^{i-1}\left(\begin{array}{c}i\\ k\end{array}\right)(i-k)\sin((i-k)\pi/2)x^{k}y^{i-k-1}z\\\noalign{\medskip}
\sum_{k=0}^i\left(\begin{array}{c}i\\ k\end{array}\right)\sin((i-k)\pi/2)x^{k}y^{i-k}
\end{array}\right).
\]
So $X_i=\lambda_1X_i^1+\lambda_2 X_i^2$ gives and explicit expression for the vector field $X_i$.

In order to investigate scenarios where the spectrum of $d^2f(p)$ is of one of the categories described in Theorem~\ref{th:2}, we have to consider more equations within the hierarchy. Following this line of reasoning, we arrive at the following result, which is our second main contribution.

\begin{theorem}[Main Theorem 2]\label{th4} Let $f$ be a smooth function with a non-degenerate critical point~$p$.  Assume that the Taylor series expansion of $f$ at $p$ has the form $f=f_0+f_2+f_4+\mathcal{O}(6)$, meaning that the homogeneous terms of order $3$ and $5$ of the function $f$ vanish. Then the unique solution to Equations~\eqref{beltrami} is $X\equiv 0$.
\end{theorem}
\noindent\textbf{Proof:}  We proceed by contradiction. Let us assume the existence of a non-trivial solution $X$ to the equations $curl(X)=fX$,
$div(X)=0$. Next, let $X_i$ be the first non-trivial term in the Taylor
series expansion of $X$ at the point $p$. As computed previously, $X_i$ satisfies the system of equations~\eqref{equation_new}, and the eigenvalues of $d^2f(p)$ are of one of the forms stated in Proposition~\ref{con2}. Let us focus on the case $i\geq 3$ and hence the eigenvalues are of the form $\{\alpha,\alpha,-i\alpha\}$. The cases $i=1$ and $i=2$ are discussed at the beginning of the proofs of Propositions~\ref{con3} and~\ref{con4}.
We distinguish  two cases: $f_0=0$ and $f_0\neq 0$. 

$\bullet$ The $\mathbf{f_0=0}$ \textbf{case:} It is easy to see that $X_i$ and $X_{i+3}$ have to satisfy the following system of equations
\begin{align*}
&  curl(X_i)=0,\\
&  div(X_{i})=0,\\
&  \langle \nabla(\sigma_1 x^2+\sigma_2 y^2+\sigma_3 z^2), X_{i} \rangle=0, \\
&  curl(X_{i+3})=(\sigma_1 x^2+\sigma_2 y^2+\sigma_3 z^2)X_i, \\
&  div(X_{i+3})=0, \\
& \langle \nabla(\sigma_1 x^2+\sigma_2 y^2+\sigma_3 z^2), X_{i+3} \rangle=0.
\end{align*}
By  Proposition~\ref{con3} in Appendix~\ref{app:con3}, we deduce that $X_i= X_{i+3}\equiv 0$, and we have completed the proof.

\vspace*{0.25cm}

$\bullet$ The $\mathbf{f_0\neq 0}$ \textbf{case:} 
We observe that $X_i$ and $X_{i+1}$ have to satisfy the system of equations
\[
\begin{array}{l}curl(X_i)=0,\\\noalign{\medskip} div(X_{i})=0,\\\noalign{\medskip}\langle \nabla(\sigma_1 x^2+\sigma_2 y^2+\sigma_3 z^3), X_{i} \rangle=0,\\\noalign{\medskip}
curl(X_{i+1})=f_0X_i,
\\\noalign{\medskip}
div(X_{i+1})=0,
\\\noalign{\medskip}
\langle \nabla(\sigma_1 x^2+\sigma_2 y^2+\sigma_3 z^2), X_{i+1} \rangle=0,
\end{array}
\]
and a computation similar to the previous one (see Proposition~\ref{con4} in Appendix~\ref{app:con4} for the details) yields the result.
\qed

\begin{remark}\label{remark} When $f_0\neq 0$ the same reasoning implies that a result analogous to Theorem~\ref{th4} holds for functions of the form
  $f=f_0+f_2+\mathcal{O}(4)$.
\end{remark}


\section{Final Remark: Perturbations}\label{conclusions}


This final section intends to outline a strategy using perturbations to strengthen our results. Since having a trivial kernel is a stable property of linear systems, we can consider perturbations of the system under consideration, which we re-write below:
\begin{align*}
&curl(X_i)=0,
  \\
  &div(X_{i})=0,
  \\
  & \langle \nabla(\sigma_1 x^2+\sigma_2 y^2+\sigma_3 z^2), X_{i} \rangle=0,
  \\
&curl(X_{i+1})=f_0X_i,
\\
&div(X_{i+1})=0,
\\
&\langle \nabla(\sigma_1 x^2+\sigma_2 y^2+\sigma_3 z^2), X_{i+1} \rangle=0.
\end{align*}
For instance, when $f = f_0 + f_2 + f_3 + \mathcal{O}(4)$ and $f_0\neq 0$ we can consider the perturbation 
\begin{align*}
&curl(X_i)=0,
  \\
  &div(X_{i})=0,\\
 & \langle \nabla(\sigma_1 x^2+\sigma_2 y^2+\sigma_3 z^2), X_{i} \rangle=0,
  \\
&curl(X_{i+1})=f_0X_i,
\\
&div(X_{i+1})=0,
\\
&\langle \nabla(\sigma_1 x^2+\sigma_2 y^2+\sigma_3 z^2), X_{i+1}
  \rangle+\epsilon\langle \nabla f_3, X_{i}\rangle=0.
\end{align*}
For $\epsilon$ small enough, the only solution to this system of equations is trivial provided that this is the case for $\epsilon=0$. This allows us to prove results like the
following one.
\begin{proposition} For a function $f$ having a non-degenerate
  critical point $p$, non vanishing at $p$ and whose term $f_3$ is
  small enough (in the sense that the coefficients are close enough to zero), the only solution to the equations $curl(X)=fX$,
  $div(X)=0$ is $X\equiv 0$.
\end{proposition}

\section*{Acknowledgments}
This work has received funding from the grants CEX2019-000904-S, RED2022-134301-T, PID2022-136795NB-I00 (D.P.-S.) and PID2022-137909NB-C21 (M.V.) funded by MCIN/AEI, and Ayudas Fundaci\'on BBVA a Proyectos de Investigaci\'on Cient\'ifica 2021 (D.P.-S.).

\newpage

\begin{appendices}
\section{$f$ with Non-isolated Critical Points}\label{nonisolated}

In this appendix we provide an example of a factor $f$ having a curve of critical points, while still admitting a non-trivial solution to Equations~\eqref{beltrami}. Namely, the function in question is
\begin{equation}\label{eqappa}
    f(x,y,z) = x^2 + y^2,
\end{equation}
where the $z$-axis constitutes a family of critical points. It is convenient to use cylindrical coordinates to construct solutions to Equations~\eqref{beltrami}. Therefore, if \[X=X^r(r,\varphi,z)e_r+X^\varphi(r,\varphi,z)e_\varphi+X^z(r,\varphi,z)e_z,\] with $\{e_r,e_\varphi,e_z\}$ the unitary cylindrical basis, then $ X^r = 0$ as a consequence of $\langle \nabla f, X \rangle =0$. Equations~\eqref{beltrami} read now:
   \begin{itemize}
       \item $curl(X) = fX$ gives
\[
       \begin{array}{l}
         \partial_rX^z - \partial_z X^\varphi = 0,
         \\ \noalign{\medskip}
        -\partial_rX^z=r^2X^\varphi,
         \\ \noalign{\medskip}
        \partial_rX^\varphi + X^\varphi/r = r^2X^z.
         \\ \noalign{\medskip}

       \end{array}
     \]
       \item $div(X) = 0$ yields
       \[
       \partial_\theta X^\theta/r + \partial_zX^z  = 0.
       \]
   \end{itemize}
One may check that the following expression provides a solution to~\eqref{beltrami}:
       \begin{flalign*}
X^r &= 0,
        \\\noalign{\bigskip}
 X^\varphi (r,\varphi ,z)  &= \frac{r \Gamma \left(\frac{2}{3}\right) J_{\frac{2}{3}}\left(\frac{r^3}{3}\right)}{\sqrt[3]{6}},
          \\\noalign{\bigskip}
         X^z(r,\varphi ,z) &=\frac{r \Gamma \left(\frac{2}{3}\right) J_{-\frac{1}{3}}\left(\frac{r^3}{3}\right)}{\sqrt[3]{6}}.  
\end{flalign*}
Here $\Gamma$ denotes the gamma function and $J$ is the Bessel function of the first kind. Using the Taylor expansion of Bessel functions at the origin, it is straightforward to check that the vector field $X$ above is analytic. In summary, we have shown that $X$ is a Beltrami field in $\mathbb R^3$ with the proportionality factor $f$ in Equation~\eqref{eqappa}. 
\section{Proof of Proposition~\ref{con2}} \label{app:con2}
\noindent\textbf{Proof:} We treat separately the cases $i = 1$, $i=2$ and $i \geq 3$.

\vspace*{0.25cm}
\noindent - \textbf{\textit{Case $i = 1$:}} Let 
\[
X_{1}=\left(\begin{array}{c}
a^{(1,0,0)}x + a^{(0,1,0)}y + a^{(0,0,1)}z\\\noalign{\medskip}
b^{(1,0,0)}x + b^{(0,1,0)}y + b^{(0,0,1)}z\\\noalign{\medskip}
c^{(1,0,0)}x + c^{(0,1,0)}y + c^{(0,0,1)}z
\end{array}\right).
\]
Then, when the coefficients of the monomials are equated to zero, the system of Equations~\eqref{equation_new} yields the following constraints:

\begin{itemize}
\item The equation $curl(X_{1})=0$ reads
\begin{equation}\label{eq:1}
    \begin{aligned}
&     -b^{(0,0,1)} + c^{(0,1,0)} = 0, \\
& a^{(0,0,1)} - c^{(1,0,0)} = 0, \\
& -a^{(0,1,0)} + b^{(1,0,0)} = 0.  \\
\end{aligned}
\end{equation}
\item The equation $div(X_{1})=0$ reads
\[
a^{(1,0,0)} + b^{(0,1,0)} + c^{(0,0,1)} = 0.
\]
\item The equation $\langle \nabla(\sigma_1 x^2+\sigma_2 y^2+\sigma_3 z^2), X_{1} \rangle=0$ reads
\begin{equation}\label{eq:2}
\begin{aligned}
 & a^{(0,1,0)} \sigma_1+b^{(1,0,0)} \sigma_2= 0, \\
& b^{(0,0,1)} \sigma_2+c^{(0,1,0)} \sigma_3 = 0,\\
 & a^{(0,0,1)} \sigma_1+c^{(1,0,0)} \sigma_3= 0, \\
 & c^{(0,0,1)} \sigma_3 = 0, \\
 & b^{(0,1,0)} \sigma_2 = 0,\\
 & a^{(1,0,0)} \sigma_1 = 0.\\
\end{aligned}
\end{equation}
\end{itemize}
Since we assumed $\sigma_i \neq 0 $ for all $i$, then $a^{(1,0,0)} = b^{(0,1,0)} = c^{(0,0,1)} = 0$ by the last three equations in~\eqref{eq:2}. Therefore, if $X_1 \neq 0$, at least two of the coefficients involved in Equations~\eqref{eq:1} must be non-zero. For instance, assuming $a^{(0,1,0)} \neq 0$, then by using the last equation in  Equations~\eqref{eq:1}, we find $b^{(1,0,0)} = a^{(0,1,0)} \neq 0$. Substituting this into the first equation in Equations~\eqref{eq:2}, we obtain $\sigma_1 + \sigma_2 = 0 \Leftrightarrow \sigma_1 = -\sigma_2$. Similar reasoning applied to the other coefficients yields the result.

\vspace*{0.25cm}
\noindent- \textbf{\textit{Case $i = 2$:}} Let 
\[
X_{2}=\left(\begin{array}{c}
a^{(2,0,0)}x^2 + a^{(1,1,0)}xy + a^{(0,2,0)}y^2 + a^{(0,1,1)}yz + a^{(0,0,2)}z^2 + a^{(1,0,1)}xz\\\noalign{\medskip}
b^{(2,0,0)}x^2 + b^{(1,1,0)}xy + b^{(0,2,0)}y^2 + b^{(0,1,1)}yz + b^{(0,0,2)}z^2 + b^{(1,0,1)}xz\\\noalign{\medskip}
c^{(2,0,0)}x^2 + c^{(1,1,0)}xy + c^{(0,2,0)}y^2 + c^{(0,1,1)}yz + c^{(0,0,2)}z^2 + c^{(1,0,1)}xz
\end{array}\right).
\]
Then, the system of Equations~\eqref{equation_new} yields the following constraints: 
\begin{itemize}

\item The equation $curl(X_{2})=0$ reads

\begin{equation}\label{eq:3}
    \begin{aligned}
 &c^{(0,1,1)}-2 b^{(0,0,2)} = 0,\\
 &2 c^{(0,2,0)}-b^{(0,1,1)} = 0,\\
 &c^{(1,1,0)}-b^{(1,0,1)} = 0,\\
 &a^{(0,1,1)}-c^{(1,1,0)} = 0,\\
 &a^{(1,0,1)}-2 c^{(2,0,0)} = 0,\\
 &b^{(1,0,1)}-a^{(0,1,1)} = 0,\\
 &2 b^{(2,0,0)}-a^{(1,1,0)} = 0,\\
  &2 a^{(0,0,2)}-c^{(1,0,1)} = 0,\\
  &b^{(1,1,0)}-2 a^{(0,2,0)} = 0.\\
\end{aligned}
\end{equation}
 \item The equation $div(X_{2})=0$ reads
 \begin{equation}\label{eq:4}
    \begin{aligned}
 &a^{(1,0,1)}+b^{(0,1,1)}+2 c^{(0,0,2)} = 0,\\
 &a^{(1,1,0)}+2 b^{(0,2,0)}+c^{(0,1,1)} = 0,\\
 &2 a^{(2,0,0)}+b^{(1,1,0)}+c^{(1,0,1)} = 0.\\
\end{aligned}
\end{equation}
\item The equation $\langle \nabla(\sigma_1 x^2+\sigma_2 y^2+\sigma_3 z^2), X_{2} \rangle=0$ reads
\begin{equation}\label{eq:5}
\begin{aligned}
& a^{(0,0,2)} \sigma_1+c^{(1,0,1)} \sigma_3 = 0,\\
& a^{(0,2,0)} \sigma_1+b^{(1,1,0)} \sigma_2 = 0,\\
& a^{(0,1,1)} \sigma_1+b^{(1,0,1)} \sigma_2+c^{(1,1,0)} \sigma_3 = 0,\\
& a^{(1,0,1)} \sigma_1+c^{(2,0,0)} \sigma_3 = 0,\\
& a^{(1,1,0)} \sigma_1+b^{(2,0,0)} \sigma_2 = 0,\\
&c^{(0,0,2)} \sigma_3 = 0,\\
& b^{(0,2,0)} \sigma_2 = 0,\\
& a^{(2,0,0)} \sigma_1 = 0.\\
\end{aligned}
\end{equation}
\end{itemize}
First, note that since $\sigma_i \neq 0 $ for all $i$, it follows that $a^{(2,0,0)} = b^{(0,2,0)} = c^{(0,0,2)} = 0$ as implied by the last three equations in~\eqref{eq:5}. Next, the key observation is the existence of a partition within the set of remaining coefficients. The subsets of this partition can be studied independently, and are given by
\begin{align*}
  &  \{a^{(0,2,0)},a^{(0,0,2)},b^{(1,1,0)},c^{(1,0,1)}\}, \
    \{a^{(1,1,0)},b^{(2,0,0)},b^{(0,0,2)},c^{(0,1,1)}\},\\
  &  \{a^{(0,1,1)},b^{(1,0,1)},c^{(1,1,0)}\},\ \text{and } 
    \{a^{(1,0,1)},b^{(0,1,1)},c^{(2,0,0)},c^{(0,2,0)}\}.
\end{align*}
We will work out the details for the first set, $\{a^{(0,2,0)},a^{(0,0,2)},b^{(1,1,0)},c^{(1,0,1)}\}$, with the other sets being analogous. Let us assume $a^{(0,2,0)} \neq 0$. Then, the last equation in~\eqref{eq:3}  shows $b^{(1,1,0)} = 2a^{(0,2,0)}$. Using this information, along with $a^{(2,0,0)} = 0$, the last equation in~\eqref{eq:4} gives $c^{(1,0,1)} = -b^{(1,1,0)} = -2a^{(0,2,0)}$. 

Next, by employing the penultimate equation in~\eqref{eq:3}, we obtain $a^{(0,0,2)} =c^{(1,0,1)}/2$. In this way, if any of the coefficients in the set  is different from zero, all of them are. Substituting into the first two equations in~\eqref{eq:5}, we get:
\begin{align*}
  c^{(1,0,1)}/2\sigma_1 +   c^{(1,0,1)}\sigma_3 = 0 \Leftarrow -\sigma_1/2 = \sigma_3, \\
  a^{(0,2,0)}\sigma_1 + 2 a^{(0,2,0)}\sigma_2 = 0\Leftrightarrow -\sigma_1/2 = \sigma_2.
\end{align*}
Therefore, $\sigma_1 + \sigma_2 + \sigma_3 = \sigma_1 -\sigma_1/2 -\sigma_1/2 = 0$. A similar computation yields the result for the remaining sets of coefficients.

\vspace*{0.25cm}
\noindent- \textbf{\textit{Case $i \geq 3$:}}  Since $\sigma_1$, $\sigma_2$ and $\sigma_3$ have different signs, otherwise by Proposition~\ref{samesign} the only solution to Equations~\eqref{equation_new} is the trivial one, we can assume that $\sigma_3$ is negative and $\sigma_1,\ \sigma_2$ are positive. If this is not the case, a change of coordinates and multiplication by $-1$ in the last equation in~\eqref{equation_new} will adjust the equations accordingly. Furthermore, based on this reasoning, we can assume that $\sigma_3=-1$, as $\sigma_1$, $\sigma_2$ and $\sigma_3$ are only determined up to a multiplicative factor.
 
 Since $curl(X_i)=0$, then $X_i=\nabla g_i$ where $g_i$ is a homogeneous polynomial of degree $i+1$. The expression $div(X_{i})=0$ now reads as $div(\nabla g_i) = \Delta g_i=0$, and so $g_i$ is a harmonic polynomial.  Then, the proposition we wish to prove can be reformulated as follows: If $g_i$ is a homogeneous harmonic polynomial of degree $i+1$ and also a first integral of a vector field of the form $(\sigma_1 x,\  \sigma_2 y,\ \sigma_3 z)$, then the set $\{\sigma_1,\sigma_2,\sigma_3\}$ must be equal to $\{ \alpha, \alpha, -i\alpha\}$. 
 
 We divide the proof  into two steps. In the first step
we demonstrate that $g_i$ has a very
particular form, $g_i=p(x,y)\cdot{}z$, where $p(x,y)$ is a harmonic polynomial on the plane. In the second step, we utilize the
expression for $X_i$ obtained in the first step to conclude the desired result. Namely, determining the values $\sigma_1, \ \sigma_2$ and $\sigma_3$ up to a multiplicative factor.

\vspace*{0.25cm}
\noindent\textbf{First step:} We  prove that $g_i$ must be of the form $p(x,y)z$ where $p(x,y)$ is a homogeneous harmonic polynomial of degree $i$ in the variables $(x,y)$. Assume that $g_i$ has the form
\[
g_i(x,y,z)=\sum_{k_1+k_2+k_3=i+1}g_i^{(k_1,k_2,k_3)}x^{k_1}y^{k_2}z^{k_3},
\]
where we will omit the subscript in $\sum$ when obvious. We impose the equations:

$\bullet$ {\it Harmonicity:} $\ div(X_i)=0\Leftrightarrow \Delta g_i =0  \Leftrightarrow\Delta \left(\sum g_i^{(k_1,k_2,k_3)}x^{k_1}y^{k_2}z^{k_3}\right)=0$ or
\begin{equation}\label{dos1} \tag{Harm}
\begin{aligned}
&\sum_{k_1 + k_2 + k_3 = i+1, k_1 \geq 2} k_1(k_1-1)g_i^{(k_1,k_2,k_3)}x^{k_1-2}y^{k_2}z^{k_3}
\\
&+\sum_{k_1 + k_2 + k_3 = i+1, k_2 \geq 2} k_2(k_2-1)g_i^{(k_1,k_2,k_3)}x^{k_1}y^{k_2-2}z^{k_3}
\\ 
&+\sum_{k_1 + k_2 + k_3 = i+1, k_3 \geq 2} k_3(k_3-1)g_i^{(k_1,k_2,k_3)}x^{k_1}y^{k_2}z^{k_3-2}=0.
\end{aligned}
\end{equation}

$\bullet$ {\it First integral:} $\
\langle\left(\begin{array}{c}\sigma_1 x\\ \sigma_2 y\\\sigma_3
               z\end{array}\right),X_i \rangle= 0\Leftrightarrow  \langle\left(\begin{array}{c}\sigma_1 x\\\sigma_2 y\\- z\end{array}\right),\nabla g_i \rangle=0$ or
\begin{equation}\label{uno1} \tag{FI}
\sum_{k_1+k_2+k_3=i+1}(\sigma_1 k_1+\sigma_2 k_2-k_3)g_i^{(k_1,k_2,k_3)}x^{k_1}y^{k_2}z^{k_3}.
\end{equation}
In what follows, we show that all the coefficients of $g_i$
that do not take the form $g^{(k_1,k_2,1)}$ vanish. This, obviously, will establish the first step of our proof. 

\vspace*{0.2cm}
{\it - Coefficients of the form $\mathbf{g^{(k_1,k_2,0)} vanish:}$ } Direct inspection of Equation~\eqref{uno1}
reveals that the terms of the form $g_i^{(k_1,k_2,0)}$ (where $k_1+k_2=i+1$) must vanish. More precisely, when all monomials in~\eqref{uno1} are equated to zero, the coefficients of the monomials of the form $x^{k_1}y^{k_2}z^{0}$ give
\[
(\sigma_1 k_1+\sigma_2 k_2)g_i^{(k_1,k_2,0)} = 0.
\]
Since $(\sigma_1 k_1+\sigma_2 k_2)$ is always positive (remember that $\sigma_1$ and  $\sigma_2$ where chosen to be positive and $k_1, \ k_2\in\mathbb{N}$, $k_1 + k_2 = i + 1\geq 4$), then $(\sigma_1 k_1+\sigma_2 k_2)g_i^{(k_1,k_2,0)}$ can only vanish if $g_i^{(k_1,k_2,0)}=0$.
\vspace*{0.2cm}

{\it -Coefficients of the form $\mathbf{g^{(k_1,k_2,2)}}$ vanish:} The
coefficients of the monomials of the form $x^{k_1}y^{k_2}z^0$ in~\eqref{dos1} are easily seen to be 
\[
(k_1+2)(k_1+2-1)g_i^{(k_1+2,k_2,0)}
+(k_2+2)(k_2+2-1)g_i^{(k_1,k_2+2,0)}
+2g_i^{(k_1,k_2,2)}.
\]
Since these coefficients must vanish by~\eqref{dos1},  and as established in the previous paragraph $g_i^{(k_1+2,k_2,0)}=g_i^{(k_1,k_2+2,0)}=0$,  it follows that $g_i^{(k_1,k_2,2)}$ also equals $0$.

\vspace*{0.2cm}

{\it -Coefficients of the form $\mathbf{g^{(k_1,k_2,3)}}$ vanish:} The
coefficients of the monomials of the form $x^{k_1}y^{k_2}z^1$ in
\eqref{dos1} can be easily identified as
\[
(k_1+2)(k_1+2-1)g_i^{(k_1+2,k_2,1)}
+(k_2+2)(k_2+2-1)g_i^{(k_1,k_2+2,1)}
+6 g_i^{(k_1,k_2,3)}.
\]
Let us assume that $g_i^{(k_1,k_2,3)}\neq 0$. Then either $g_i^{(k_1+2,k_2,1)}$ or $g_i^{(k_1,k_2+2,1)}$ must be different from zero, as the last expression must vanish due to~\eqref{dos1}. Suppose that $g_i^{(k_1+2,k_2,1)}\neq 0$ and the other case can be treated in an analogous way.  Now, the coefficients of the monomials $x^{k_1}y^{k_2}z^3$ and $x^{k_1+2}y^{k_2}z$ in~\eqref{uno1} are $(\sigma_1k_1+\sigma_2 k_2-3)g_i^{(k_1,k_2,3)}$ and $(\sigma_1(k_1+2)+\sigma_2k_2-1)g_i^{(k_1+2,k_2,1)}$ respectively. Since by Equation~\eqref{uno1} these coefficients must vanish, the only possibility is
\begin{align*}
&\sigma_1k_1+\sigma_2 k_2-3=0, \\
&\sigma_1(k_1+2)+\sigma_2k_2-1=0.
\end{align*}
The second equation minus the first one gives 
$2\sigma_1+2=0\Leftrightarrow \sigma_1=-1$,
which is a contradiction as we assumed that $\sigma_1$ is
positive. Therefore, we conclude that the coefficients of the form
$g^{(k_1,k_2,3)}$ must vanish. Notice that we are using the fact that $k_1$ and $k_2$ cannot vanish at the same time, as $k_1 + k_2 + k_3 = i + 1$, and since $i\geq 3$ and $k_3 = 3$ then $k_1 + k_2 \geq 1$.

\vspace*{0.2cm}

{\it -Coefficients of the form $\mathbf{g^{(k_1,k_2,j)}}$ with $j>
  3$ vanish:} We proceed by induction. Let us assume that we have
proved the coefficients of the form $g_i^{(k_1,k_2,k_3)}$ vanish for
$2\leq k_3\leq j-1$. Using  Equation~\eqref{dos1} the coefficients of the form $g_i^{(k_1,k_2,j)}$ must satisfy the equation
\begin{align*}
 &(k_1+2)(k_1+2-1)g_i^{(k_1+2,k_2,j-2)}
+(k_2+2)(k_2+2-1)g_i^{(k_1,k_2+2,j-2)}
\\
&+j(j-1) g_i^{(k_1,k_2,j)}=0.   
\end{align*}
 By the induction hypothesis,  $g_i^{(k_1+2,k_2,j-2)}=g_i^{(k_1,k_2+2,j-2)}=0$,
which implies $g_i^{(k_1,k_2,j)}=0$. Since the only coefficients that do not
vanish are of the form $g_i^{(k_1,k_2,1)}$,  the first step is
complete. The fact that $p(x,y)$ is harmonic is obvious.

\vspace*{0.2cm}

\noindent\textbf{Second step:}  In this step we obtain a
linear system of equations involving $\sigma_1$ and $\sigma_2$, which determines their  values once the value of $\sigma_3$
is fixed. More precisely, we obtain $\sigma_1=\sigma_2=1/i$ when $\sigma_3=-1$. Remember $g_i=p(x,y)\cdot{}z$ where $p(x,y)$ is a harmonic polynomial of degree $i$ in two variables. Thus, using the same notation as before
\[
p(x,y)=\sum_{k_1+k_2=i}p^{(k_1,k_2)}x^{k_1}y^{k_2}.
\]
Since $g_i$  satisfies the
Equation~\eqref{uno1}, then
\[
\langle (\sigma_1 x,\sigma_2 y,- z),\nabla g_i \rangle=\sum(\sigma_1 k_1+\sigma_2 k_2-1)p^{(k_1,k_2)}x^{k_1}y^{k_2}z=0,
\]
where we can disregard here the variable $z$ and we obtain
\begin{equation}\label{tres}
\sum(\sigma_1 k_1+\sigma_2 k_2-1)p^{(k_1,k_2)}x^{k_1}y^{k_2}=0.
\end{equation}
Considering that $i\geq 3$, in all monomials of
$p(x,y)$ either $k_1$ or $k_2$ must be greater or equal than $2$. Given that we are assuming $p(x,y)\neq 0$, there exist integers
$\hat{k}_1$, $\hat{k}_2$ such that $p^{(\hat{k}_1,\hat{k}_2)}\neq 0$. Let us
assume that $\hat{k}_1\geq 2$, the other case is treated
analogously. Using that
\[
\Delta p(x,y)=0\Leftrightarrow \Delta\left( \sum_{k_1+k_2=i}p^{(k_1,k_2)}x^{k_1}y^{k_2}\right)= 0,
\]
we have
\[
\sum k_1(k_1-1)p^{(k_1,k_2)}x^{k_1-2}y^{k_2}
+\sum k_2(k_2-1)p^{(k_1,k_2)}x^{k_1}y^{k_2-2}=0.
\]
In this equation the coefficient of the monomial $x^{\hat{k}_1-2}y^{\hat{k}_2}$ is 
\[
\hat{k}_1(\hat{k}_1-1)p^{(\hat{k}_1,\hat{k}_2)}+(\hat{k}_2+2)(\hat{k}_2+2-1)p^{(\hat{k}_1-2,\hat{k}_2+2)},
\]
which, by $\Delta p(x,y)=0$, has to vanish. Therefore, if $p^{(\hat{k}_1,\hat{k}_2)}\neq 0$ then $p^{(\hat{k}_1-2,\hat{k}_2+2)}\neq 0$. Going back to Equations~\eqref{tres}, and equating to zero the coefficients of the monomials  $x^{\hat{k}_1}y^{\hat{k}_2}$ and $x^{\hat{k}_1-2}y^{\hat{k}_2+2}$, we obtain
$
(\hat{k}_1\sigma_1+\hat{k}_2\sigma_2-1)p^{(\hat{k}_1,\hat{k}_2)} = 0$ and 
$((\hat{k}_1-2)\sigma_1+(\hat{k}_2+2)\sigma_2-1)p^{(\hat{k}_1-2,\hat{k}_2+2)} = 0$. Since $p^{(\hat{k}_1,\hat{k}_2)}\neq 0$ and $p^{(\hat{k}_1-2,\hat{k}_2+2)}\neq 0$, $\sigma_1$ and $\sigma_2$ have to satisfy the following  system of linear equations
\begin{align*}
& \hat{k}_1\sigma_1+\hat{k}_2\sigma_2=1,\\ 
&(\hat{k}_1-2)\sigma_1+(\hat{k}_2+2)\sigma_2=1.
\end{align*}
As the determinant of this system is
$
\hat{k}_1(\hat{k}_2+2)-(\hat{k}_1-2)\hat{k}_2=\hat{k}_1\hat{k}_2+2\hat{k}_1-\hat{k}_1\hat{k}_2+2\hat{k}_2=2(\hat{k}_1+\hat{k}_2)=2i\neq 0,
$ there is a unique solution to the last system of equations. Due to
$\hat{k}_1+\hat{k}_2=i$, it is easy to see that $\sigma_1=\sigma_2=1/i$ is the
only solution. Summarizing, when $\sigma_3=-1$ then
$\sigma_1=\sigma_2=1/i$. Since solutions for $\sigma_1$, $\sigma_2$ and $\sigma_3$ are obtained up to a multiplicative
factor, the desired result follows.
\qed

\section{Proposition~\ref{con3}}\label{app:con3} 
The following result is instrumental in the proof of Theorem~\ref{th4}:

\begin{proposition}\label{con3} Let $X_i$, $X_{i+3}$ be  vector fields whose components are
  homogeneous polynomials of degree $i$ and $i +3$ respectively. Then, the non-trivial solutions $X_i$, $X_{i + 3}$ to
  the equations
  \begin{subequations}\label{aaa}
    \begin{align}
&  curl(X_i)=0, \label{need4}\\ 
&  div(X_{i})=0,\\
  &  \langle \nabla(\sigma_1 x^2+\sigma_2 y^2+\sigma_3 z^2), X_{i} \rangle=0, \label{need5}
  \\
  &  curl(X_{i+3})=(\sigma_1 x^2+\sigma_2 y^2+\sigma_3 z^2)X_i, \label{need1}
  \\
&  div(X_{i+3})=0, \label{need2}
\\
      & \langle \nabla(\sigma_1 x^2+\sigma_2 y^2+\sigma_3 z^2), X_{i+3} \rangle=0, \label{need3}
    \end{align}
\end{subequations}
with non-zero $\sigma_1,\ \sigma_2, \ \sigma_3$ (having different signs), are of the form $X_i\equiv 0$ and (up to a permutation of the variables) $X_{i+3}=\nabla(p(x,y)z)$, where $p(x,y)$ is a harmonic polynomial of degree $i + 3$ on the plane. Moreover, non-trivial solutions only exist if $\{\sigma_1,\ \sigma_2, \ \sigma_3\}=\{\alpha,\alpha,-(i+3)\alpha\}$ where $\alpha$ is a non-zero real number.
\end{proposition}

\noindent\textbf{Proof:} The cases $i=1,\ 2$ follow by a lengthy but straightforward
computation. For $i = 1, \ 2$ we also include the code for the corresponding symbolic computations using Mathematica, which can be accessed via this  \href{https://www.dropbox.com/scl/fo/fzj8wckeucanp47k5fhfq/h?rlkey=t25vlwf7bfv3dgf78elt60sy6&dl=0}{link}. Therefore, we shall focus on the case $i\geq 3$. 

Let us assume that $X_i$ is non-zero, a direct application of Proposition~\ref{con2} to the sub-system~\eqref{need4}--\eqref{need5} yields  $X_i=\lambda_1X_i^1+\lambda_2X_i^2$ and $\{\sigma_1,\sigma_2,\sigma_3\} = \{\alpha,\alpha,-i\alpha\}$. Moreover, we can assume that
$\sigma_1=\sigma_2=\alpha$ and $\sigma_3=-i\alpha$ for some positive real number $\alpha$, otherwise we can make
a change of coordinates to get to this situation. Furthermore, we can take $\alpha = 1$ by simply multiplying by $\alpha^{-1}$ in~\eqref{need1} and~\eqref{need3}.

We proceed by contradiction. Our strategy relies on demonstrating that if $X_i\not\equiv 0$, then the system~\eqref{aaa} is incompatible. Once we show $X_i \equiv 0$, application of Proposition~\ref{con2} to the Equations~\eqref{need1}--\eqref{need3} gives the desired result for the possible values of  the $\sigma$'s and $X_{i + 3}$. We introduce the notation:
\[
X_{i+3}=\left(\begin{array}{c}\sum_{k_1+k_2+k_3=i+3}a^{(k_1,k_2,k_3)}x^{k_1}y^{k_2}z^{k_3}\\\noalign{\medskip}
\sum_{k_1+k_2+k_3=i+3}b^{(k_1,k_2,k_3)}x^{k_1}y^{k_2}z^{k_3}\\\noalign{\medskip}
\sum_{k_1+k_2+k_3=i+3}c^{(k_1,k_2,k_3)}x^{k_1}y^{k_2}z^{k_3}
\end{array}\right),
\]
and Equations~\eqref{need1},~\eqref{need2} and~\eqref{need3} read:

\vspace*{0.2cm}


$\bullet$ {\it Rotational:} $curl(X_{i+3})=(\sigma_1 x^2+\sigma_2
y^2+\sigma_3 z^2)X_i\Leftrightarrow curl(X_{i+3})=(\sigma_1 x^2+\sigma_2
y^2+\sigma_3 z^2)\lambda_1X_i^1+(\sigma_1 x^2+\sigma_2
y^2+\sigma_3 z^2)\lambda_2X_i^2$ or
\begin{equation}\label{rotational} \tag{Rot-P7}
\begin{array}{l}
\left(
\begin{array}{c}
(\sum k_2c^{(k_1,k_2,k_3)}x^{k_1}y^{k_2-1}z^{k_3}-\sum k_3b^{(k_1,k_2,k_3)})x^{k_1}y^{k_2}z^{k_3-1}\\\noalign{\medskip}
(\sum k_3a^{(k_1,k_2,k_3)}x^{k_1}y^{k_2}z^{k_3-1}-\sum k_1c^{(k_1,k_2,k_3)})x^{k_1-1}y^{k_2}z^{k_3}\\\noalign{\medskip}
(\sum k_1b^{(k_1,k_2,k_3)}x^{k_1-1}y^{k_2}z^{k_3}-\sum k_2a^{(k_1,k_2,k_3)})x^{k_1}y^{k_2-1}z^{k_3}
\end{array}\right)
\\ \noalign{\bigskip}
=
(\sigma_1 x^2+\sigma_2 y^2+\sigma_3 z^2)\lambda_1\left(\begin{array}{c}\sum_{k=1}^i\left(\begin{array}{c}i\\ k\end{array}\right)k\cos((i-k)\pi/2)x^{k-1}y^{i-k}z \\\noalign{\medskip}
\sum_{k=0}^{i-1}\left(\begin{array}{c}i\\ k\end{array}\right)(i-k)\cos((i-k)\pi/2)x^{k}y^{i-k-1}z\\\noalign{\medskip}
\sum_{k=0}^i\left(\begin{array}{c}i\\ k\end{array}\right)\cos((i-k)\pi/2)x^{k}y^{i-k}
\end{array}\right)
\\ \noalign{\bigskip}
+(\sigma_1 x^2+\sigma_2 y^2+\sigma_3 z^2)\lambda_2\left(\begin{array}{c}\sum_{k=1}^i\left(\begin{array}{c}i\\ k\end{array}\right)k\sin((i-k)\pi/2)x^{k-1}y^{i-k}z \\\noalign{\medskip}
\sum_{k=0}^{i-1}\left(\begin{array}{c}i\\ k\end{array}\right)(i-k)\sin((i-k)\pi/2)x^{k}y^{i-k-1}z\\\noalign{\medskip}
\sum_{k=0}^i\left(\begin{array}{c}i\\ k\end{array}\right)\sin((i-k)\pi/2)x^{k}y^{i-k}
\end{array}\right).
\end{array}
\end{equation}

\vspace*{0.2cm}
$\bullet$ {\it Divergence:} $div(X_{i+3})=0$ or
\begin{equation}\label{divergence} \tag{Div-P7}
\begin{aligned}
&\sum_{k_1+ k_2 + k_3 = i + 3,k_1\geq 1} k_1a^{(k_1,k_2,k_3)}x^{k_1-1}y^{k_2}z^{k_3}\\
&+ \sum_{k_1+ k_2 + k_3 = i + 3,k_2\geq 2}k_2b^{(k_1,k_2,k_3)}x^{k_1}y^{k_2-1}z^{k_3}\\
&+\sum_{k_1+ k_2 + k_3 = i + 3,k_3\geq 1}k_3c^{(k_1,k_2,k_3)}x^{k_1}y^{k_2}z^{k_3-1}=0    .
\end{aligned}
\end{equation}

\vspace*{0.2cm}
$\bullet$ {\it First integral:} $\langle \nabla(\sigma_1 x^2+\sigma_2 y^2+\sigma_3 z^2), X_{i+3} \rangle=0$ or
\begin{equation}\label{integral} \tag{FI-P7}
\begin{aligned}
&\langle \left(\begin{array}{c}1\\\noalign{\medskip}1\\\noalign{\medskip}-i\end{array}\right),\left(\begin{array}{c}\sum_{k_1+k_2+k_3=i+3} a^{(k_1,k_2,k_3)}x^{k_1}y^{k_2}z^{k_3}\\\noalign{\medskip}
\sum_{k_1+k_2+k_3=i+3}b^{(k_1,k_2,k_3)}x^{k_1}y^{k_2}z^{k_3}\\\noalign{\medskip}
\sum_{k_1+k_2+k_3=i+3}c^{(k_1,k_2,k_3)}x^{k_1}y^{k_2}z^{k_3}
\end{array}\right) \rangle
\\ \noalign{\medskip}
&=\sum_{k_1 + k_2 + k_3 = i + 3} a^{(k_1,k_2,k_3)}x^{k_1+1}y^{k_2}z^{k_3}+ b^{(k_1,k_2,k_3)}x^{k_1}y^{k_2+1}z^{k_3}
\\
&\quad -i c^{(k_1,k_2,k_3)}x^{k_1}y^{k_2}z^{k_3+1}.
\end{aligned}
\end{equation}
From the Equation~\eqref{divergence}, taking the coefficients of the monomials $y^{i+2}$ and $xy^{i+1}$, we get the expressions
\begin{equation}\label{divergencias}
\begin{aligned}
&a^{(1,i+2,0)}+(i+3)b^{(0,i+3,0)}+c^{(0,i+2,1)}&=0,\\
&2a^{(2,i+1,0)}+(i+2)b^{(1,i+2,0)}+c^{(1,i+1,1)}&=0.
\end{aligned}
\end{equation}
Next, we apply the following two-steps strategy. First, we determine the values of the variables in Equations~\eqref{divergencias}
 as functions of $i,\ \lambda_1$ and $\lambda_2$. This task relies on solving several
subsystems of equations derived from the Equations~\eqref{rotational} and~\eqref{divergence}. Secondly, we show that the values obtained for the variables  are incompatible with Equations~\eqref{divergencias} above for all values of $\lambda_1\neq 0, \lambda_2 \neq 0$, thereby yielding the desired contradiction.

\vspace*{0.2cm}
\noindent{\it-Computation of $\mathbf{a^{(1,i+2,0)}}$:} The
coefficients of the monomials $xy^{i+1}$ and $x^2y^{i+2}$ in Equations~\eqref{rotational} and~\eqref{integral} yield respectively:
\begin{align*}
&2b^{(2,i+1,0)} - (i+2)a^{(1,i+2,0)}=i\left(\lambda_1\cos(\frac{(i-1)\pi}{2})+\lambda_2\sin(\frac{(i-1)\pi}{2})\right),\\
&a^{(1,i+2,0)}+b^{(2,i+1,0)}=0.
\end{align*}
The determinant of this system is $-i-4$ which only vanishes for $i=-4$. Since we have assumed $i\geq 3$, the solution of the system is uniquely determined and we obtain
\begin{align*}
 a^{(1,i+2,0)}=\displaystyle\frac{-i\left(\lambda_1  \cos \left(\frac{(i-1)\pi}{2}\right)+\lambda_2  \sin
   \left(\frac{(i-1)\pi}{2}\right)\right)}{i+4}.
\end{align*}

\vspace*{0.2cm}
\noindent{\it-Computation of $\mathbf{b^{(0,i+3,0)}}$:} It is easy to see from~\eqref{integral} that this coefficient must vanish.

\vspace*{0.2cm}
\noindent{\it -Computation of $\mathbf{c^{(0,i+2,1)}}$:} The coefficients of the monomials $y^{i+1}z$ and $y^{i+2}z^2$ in Equations~\eqref{rotational} and~\eqref{integral} yield respectively:
\begin{align*}
&(i+2)c^{(0,i+2,1)}-2b^{(0,i+1,2)}=i\left(\lambda_1\cos(\frac{(i-1)\pi}{2})+\lambda_2\sin(\frac{(i-1)\pi}{2})\right),\\
&b^{(0,i+1,2)}-ic^{(0,i+2,1)}=0.
\end{align*}
The determinant of this system is $i-2$. Therefore, the unique solution for $c^{(0,i+2,1)}$ is given by
\begin{align*}
&c^{(0,i+2,1)}=\displaystyle\frac{-i \left(\lambda_1 \cos \left(\frac{(i-1)\pi}{2}\right)+ \lambda_2 \sin \left(\frac{(i-1)\pi}{2}\right)\right)}{i-2}.
\end{align*}

\vspace*{0.2cm}

\noindent{\it -Computation of $\mathbf{a^{(2,i+1,0)}}$:} Taking the
coefficients of the monomials $x^2y^{i}$  in Equations~\eqref{rotational} and the monomial $x^3y^{i+1}$ in~\eqref{integral}, we obtain the system:
\begin{align*}
&3b^{(3,i,0)} -(i+1)a^{(2,i+1,0)}=(1-\left(\begin{array}{c} i \\2 \end{array}\right))(\lambda_1\cos(\frac{i\pi}{2})+\lambda_2\sin(\frac{i\pi}{2})),\\
&a^{(2,i+1,0)}+b^{(3,i,0)}=0,
\end{align*}
which has determinant $-i-4$. The system is completely determined and its solution for $a^{(2,i+1,0)}$ is given by
\[
a^{(2,i+1,0)}=\displaystyle\frac{\left(i^2-i-2\right) \left(\lambda_1 \cos \left(\frac{\pi  i}{2}\right)+\lambda_2 \sin \left(\frac{\pi
    i}{2}\right)\right)}{2 (i+4)}.
\]
\vspace*{0.2cm}

\noindent{\it -Computation of $\mathbf{b^{(1,i+2,0)}}$:} Taking the
coefficients of the  monomials $y^{i+2}$  in Equations~\eqref{rotational} and the monomial $xy^{i+3}$ in~\eqref{integral} yields the system:
\begin{align*}
&b^{(1,i+2,0)} - (i+3)a^{(0,i+3,0)}=\lambda_1\cos(\frac{i\pi}{2})+\lambda_2\sin(\frac{i\pi}{2}),\\\noalign{\medskip}
&a^{(0,i+3,0)}+b^{(1,i+2,0)}=0,
\end{align*}
which has determinant $-i-4$. Therefore, the system above is completely determined and its unique solution for $b^{(1,i+2,0)}$ is given by
\[
b^{(1,i+2,0)}=\displaystyle\frac{\lambda_1 \cos \left(\frac{i\pi }{2}\right)+\lambda_2 \sin \left(\frac{i\pi }{2}\right)}{i+4}.
\]
\vspace*{0.2cm}

\noindent{\it -Computation of $\mathbf{c^{(1,i+1,1)}}$:} Taking the monomials $y^{i+1}z$ and $xy^iz$ in Equations~\eqref{rotational} and the monomial $xy^{i+2}z$ in~\eqref{integral} we obtain the system:
\begin{align*}
&2a^{(0,i+1,2)}-c^{(1,i+1,1)}=i\left(\lambda_1\cos(\frac{i\pi}{2})+\lambda_2\sin(\frac{i\pi}{2})\right),\\
&(i+1)c^{(1,i+1,1)} -2b^{(1,i,2)} =\left(\begin{array}{c}i\\2\end{array}\right)2\left(\lambda_1\cos(\frac{(i-2)\pi}{2})+\lambda_2\sin(\frac{(i-2)\pi}{2})\right),\\
&a^{(0,i+1,2)}+b^{(1,i,2)}-ic^{(1,i+1,1)}=0,
\end{align*}
which has determinant $2i - 4$ and whose solution for $c^{(1,i+1,1)}$ is
\[
c^{(1,i+1,1)}=\displaystyle\frac{i^2 \left(\lambda_1 \cos \left(\frac{i\pi  }{2}\right)+\lambda_2 \sin \left(\frac{i\pi }{2}\right)\right)}{i-2}.
\]

The computations above allow us to determine the values of all variables involved in Equation~\eqref{divergencias}, concluding the first step of our approach. By substituting these values into~\eqref{divergencias} we get the following constraints on $\lambda_1$ and $\lambda_2$.
\small{
\begin{align*}
& \displaystyle\frac{-i\lambda_1  \cos \left(\frac{(i-1)\pi}{2}\right)-i\lambda_2  \sin
   \left(\frac{(i-1)\pi}{2}\right)}{i+4}+\displaystyle\frac{-i \lambda_1 \cos \left(\frac{(i-1)\pi}{2}\right)-i \lambda_2 \sin
   \left(\frac{(i-1)\pi}{2}\right)}{i-2}  = 0,
	\\ 
& 2\displaystyle\frac{\left(i^2-i-2\right) \left(\lambda_1 \cos \left(\frac{\pi  i}{2}\right)+\lambda_2 \sin \left(\frac{\pi
    i}{2}\right)\right)}{2 (i+4)}+(i+2)\displaystyle\frac{\lambda_1 \cos \left(\frac{\pi i}{2}\right)+\lambda_2 \sin \left(\frac{\pi i}{2}\right)}{i+4} 
    \\& +\displaystyle\frac{i^2 \left(\lambda_1 \cos \left(\frac{
\pi  i}{2}\right)+\lambda_2 \sin \left(\frac{\pi i}{2}\right)\right)}{i-2} = 0.
\end{align*}
}
\noindent After rearranging,
\begin{equation}\label{lastequ}
\begin{aligned}
&\left(\displaystyle\frac{-i}{i+4}+\displaystyle\frac{-i}{i-2}\right)\left(\lambda_1  \cos \left(\frac{(i-1)\pi}{2}\right)+\lambda_2  \sin
   \left(\frac{(i-1)\pi}{2}\right)\right)= 0,
\\ 
& \left(\displaystyle\frac{i^2-i-2 }{ i+4} +\displaystyle\frac{i+2}{i+4}+\displaystyle\frac{i^2}{i-2}\right)\left(\lambda_1 \cos \left(\frac{\pi  i}{2}\right)+\lambda_2 \sin \left(\frac{\pi
    i}{2}\right)\right)=0.
\end{aligned}
\end{equation}
Now, since for $i\geq 3$, $\left(\displaystyle\frac{-i}{i+4}+\displaystyle\frac{-i}{i-2}\right) \neq 0$ and $\left(\displaystyle\frac{i^2-i-2 }{i+4} +\displaystyle\frac{i+2}{i+4}+\displaystyle\frac{i^2}{i-2}\right) \neq 0$, Equations~\eqref{lastequ} can only vanish if
\begin{equation}\label{casi}
\begin{aligned}
&\lambda_1  \cos \left(\frac{(i-1)\pi}{2}\right)+\lambda_2  \sin
   \left(\frac{(i-1)\pi}{2}\right)=0,
	\\
&\lambda_1 \cos \left(\frac{i\pi  }{2}\right)+\lambda_2 \sin \left(\frac{i\pi
    }{2}\right)=0.
\end{aligned}
\end{equation}
Finally, we observe that depending on the parity of $i$ we have:

\vspace*{0.2cm}
{\it -For $i$ even}, the system~\eqref{casi} reads
\begin{align*}
&\lambda_2  \sin
   \left(\frac{(i-1)\pi}{2}\right)=0,
	\\
&\lambda_1 \cos \left(\frac{i\pi  }{2}\right)=0.
\end{align*}

\vspace*{0.2cm}
{\it -For $i$ odd}, the system~\eqref{casi} reads
\begin{align*}
&\lambda_1  \cos \left(\frac{(i-1)\pi}{2}\right)=0,
	\\
&\lambda_2 \sin \left(\frac{i\pi
    }{2}\right)=0.
\end{align*}
In both cases we can conclude that $\lambda_1=\lambda_2=0$, and, as a consequence, $X_i \equiv 0$. This gives the desired contradiction. A straightforward application of Proposition~\ref{con2} to the Equations~\eqref{need1}--\eqref{need3} concludes the proof.

\qed

\section{Proposition~\ref{con4}}\label{app:con4} 
The following proposition is also used in the proof of Theorem~\ref{th4}:

\begin{proposition}\label{con4}Let $X_i$, $X_{i + 1}$ be  vector fields whose components are
  homogeneous polynomials of degree $i$ and $i +1$ respectively. Then, the non-trivial solutions $X_i$, $X_{i+1}$ to
  the equations
\begin{subequations}
  \begin{align}
&curl(X_i)=0,
  \\
  &div(X_{i})=0,\\
  &\langle \nabla(\sigma_1 x^2+\sigma_2 y^2+\sigma_3 z^2), X_{i} \rangle=0,
  \\
&curl(X_{i+1})=f_0X_i, \label{eq1}
\\
&div(X_{i+1})=0,  
\\
&\langle \nabla(\sigma_1 x^2+\sigma_2 y^2+\sigma_3 z^2), X_{i+1} \rangle=0, \label{eq3}
\end{align}
\end{subequations}
with  $\sigma_1,\ \sigma_2, \ \sigma_3, \ f_0$ non-zero constants ($\sigma_1,\ \sigma_2$, and  $\sigma_3$ having different signs), are of the form $X_i\equiv 0$.
Moreover, for $i\geq 2$ non-trivial solutions only exist if $\{\sigma_1,\ \sigma_2, \ \sigma_3\}=\{\alpha,\alpha,-(i+1)\alpha\}$ where $\alpha$ is a non-zero real number.
\end{proposition}

\noindent\textbf{Proof:} We employ the same methodology outlined in Proposition~\ref{con3}. The cases $i=1,\ 2$ follow by a lengthy but straightforward
computation.  We also include the code for the corresponding symbolic computations  (cases $i = 1,\ 2$) using Mathematica, which can be accessed via this \href{https://www.dropbox.com/scl/fo/fzj8wckeucanp47k5fhfq/h?rlkey=t25vlwf7bfv3dgf78elt60sy6&dl=0}{link}. Thus, we shall focus on $i\geq 3$. We introduce the notation
\[
X_{i+1}=\left(\begin{array}{c}\sum a^{(k_1,k_2,k_3)}x^{k_1}y^{k_2}z^{k_3}\\\noalign{\medskip}
\sum b^{(k_1,k_2,k_3)}x^{k_1}y^{k_2}z^{k_3}\\\noalign{\medskip}
\sum c^{(k_1,k_2,k_3)}x^{k_1}y^{k_2}z^{k_3}
\end{array}\right).
\]
We may assume that that $\alpha=f_0=1$ without loss of generality. Moreover, we can also assume that $\sigma_1 = \sigma_2 = 1$ and $\sigma_3 = -i$. Then, the equations~\eqref{eq1}--\eqref{eq3} become:

$\bullet$ {\it Rotational:} $curl(X_{i+1})=X_i\Leftrightarrow curl(X_{i+1})=\lambda_1X_i^1+\lambda_2X_i^2$ or
\begin{equation}\label{rotational2} \tag{Rot-P8}
\begin{array}{l}
\left(
\begin{array}{c}
(\sum k_2c^{(k_1,k_2,k_3)}x^{k_1}y^{k_2-1}z^{k_3}-\sum k_3b^{(k_1,k_2,k_3)})x^{k_1}y^{k_2}z^{k_3-1}\\
(\sum k_3a^{(k_1,k_2,k_3)}x^{k_1}y^{k_2}z^{k_3-1}-\sum k_1c^{(k_1,k_2,k_3)})x^{k_1-1}y^{k_2}z^{k_3}\\
(\sum k_1b^{(k_1,k_2,k_3)}x^{k_1-1}y^{k_2}z^{k_3}-\sum k_2a^{(k_1,k_2,k_3)})x^{k_1}y^{k_2-1}z^{k_3}
\end{array}\right)
\\ \noalign{\bigskip}
=
\lambda_1\left(\begin{array}{c}\sum_{k=1}^i\left(\begin{array}{c}i\\ k\end{array}\right)k\cos((i-k)\pi/2)x^{k-1}y^{i-k}z \\\noalign{\medskip}
\sum_{k=0}^{i-1}\left(\begin{array}{c}i\\ k\end{array}\right)(i-k)\cos((i-k)\pi/2)x^{k}y^{i-k-1}z\\
\sum_{k=0}^i\left(\begin{array}{c}i\\ k\end{array}\right)\cos((i-k)\pi/2)x^{k}y^{i-k}
\end{array}\right)
\\ \noalign{\bigskip}
+\lambda_2\left(\begin{array}{c}\sum_{k=1}^i\left(\begin{array}{c}i\\ k\end{array}\right)k\sin((i-k)\pi/2)x^{k-1}y^{i-k}z \\
\sum_{k=0}^{i-1}\left(\begin{array}{c}i\\ k\end{array}\right)(i-k)\sin((i-k)\pi/2)x^{k}y^{i-k-1}z\\
\sum_{k=0}^i\left(\begin{array}{c}i\\ k\end{array}\right)\sin((i-k)\pi/2)x^{k}y^{i-k}
\end{array}\right).
\end{array}
\end{equation}

\vspace*{0.2cm}
$\bullet$ {\it Divergence:} $div(X_{i+1})=0$ or
\begin{equation}\label{divergence2} \tag{Div-P8}
\begin{aligned}
 &   \sum k_1a^{(k_1,k_2,k_3)}x^{k_1-1}y^{k_2}z^{k_3}+ \sum k_2b^{(k_1,k_2,k_3)}x^{k_1}y^{k_2-1}z^{k_3} \\ & + \sum  k_3c^{(k_1,k_2,k_3)}x^{k_1}y^{k_2}z^{k_3-1}=0.
\end{aligned}
\end{equation}

\vspace*{0.2cm}
$\bullet$ {\it First integral:} $\langle \nabla(\sigma_1 x^2+\sigma_2 y^2+\sigma_3 z^3), X_{i+1} \rangle=0$ or
\begin{equation}\label{integral2} \tag{FI-P8}
\begin{array}{c}
\langle \left(\begin{array}{c}1\\\noalign{\medskip}1\\\noalign{\medskip}-i\end{array}\right),\left(\begin{array}{c}\sum a^{(k_1,k_2,k_3)}x^{k_1}y^{k_2}z^{k_3}\\\noalign{\medskip}
\sum b^{(k_1,k_2,k_3)}x^{k_1}y^{k_2}z^{k_3}\\\noalign{\medskip}
\sum c^{(k_1,k_2,k_3)}x^{k_1}y^{k_2}z^{k_3}
\end{array}\right) \rangle
\\ \noalign{\bigskip}
=\sum a^{(k_1,k_2,k_3)}x^{k_1+1}y^{k_2}z^{k_3}+ b^{(k_1,k_2,k_3)}x^{k_1}y^{k_2+1}z^{k_3}-i c^{(k_1,k_2,k_3)}x^{k_1}y^{k_2}z^{k_3+1}.
\end{array}
\end{equation}
We apply the same two-steps strategy as in Proposition~\ref{con3}. From~\eqref{divergence2}, taking the coefficients of the monomials $y^i$ and
$xy^{i-1}$ we obtain the following system of equations:
\begin{equation}\label{divergences2}
\begin{array}{l}
a^{(1,i,0)}+(i+1)b^{(0,i+1,0)}+c^{(0,i,1)}=0,\\ \noalign{\medskip}
2a^{(2,i-1,0)}+ib^{(1,i,0)}+c^{(1,i-1,1)}=0.
\end{array}
\end{equation}
First, we compute the values of the variables in Equation~\eqref{divergences2} as functions of $i$,  $\lambda_1$ and $\lambda_2$. Secondly, we demonstrate that substituting the obtained values back into~\eqref{divergences2} leads to the conclusion that the only permissible solution is $\lambda_1 = \lambda_2 = 0$.

\vspace*{0.2cm}
\noindent{\it-Computation of $\mathbf{a^{(1,i,0)}}$:} The coefficients
of the monomials $xy^{i-1}$ and $x^2y^{i}$ in~\eqref{rotational2} and~\eqref{integral2} yield:
\begin{align*}
&2b^{(2,i-1,0)} -i a^{(1,i,0)} =i\left(\lambda_1\cos(\frac{(i-1)\pi}{2})+\lambda_2\sin(\frac{(i-1)\pi}{2})\right),\\
&a^{(1,i,0)}+b^{(2,i-1,0)}=0.
\end{align*}
The  determinant of the system above is $-i-2$ which only vanishes for $i=-2$. Since we are assuming $i\geq 3$, then the system is completely determined and the
solution for $a^{(1,i,0)}$ is given by
\begin{align*}
& a^{(1,i,0)}=\displaystyle\frac{-i\left(\lambda_1  \cos \left(\frac{(i-1)\pi}{2}\right)+\lambda_2  \sin
   \left(\frac{(i-1)\pi}{2}\right)\right)}{i+2}.
\end{align*}

\vspace*{0.2cm}
\noindent{\it-Computation of $\mathbf{b^{(0,i+1,0)}}$:} It is easy to see from~\eqref{integral2} that this coefficient has to vanish.

\vspace*{0.2cm}
\noindent{\it -Computation of $\mathbf{c^{(0,i,1)}}$:} The
coefficients of the  monomials $y^{i-1}z$ and $y^{i}z^2$ in~\eqref{rotational2} and~\eqref{integral2} yield respectively:
\begin{align*}
&  ic^{(0,i,1)}  -2b^{(0,i-1,2)}=i\left(\lambda_1\cos(\frac{(i-1)\pi}{2})+\lambda_2\sin(\frac{(i-1)\pi}{2})\right),\\
& b^{(0,i-1,2)}-ic^{(0,i,1)}=0,
\end{align*}
which has determinant $i$. Then, the systems is completely determined
 for
$i\geq 3$ and the solution for $c^{(0,i,1)}$ is given by
\begin{align*}
c^{(0,i,1)}= -\lambda_1 \cos \left(\frac{(i-1)\pi}{2}\right)- \lambda_2 \sin
   \left(\frac{(i-1)\pi}{2}\right).
\end{align*}
\vspace*{0.2cm}

\noindent{\it -Computation of $\mathbf{a^{(2,i-1,0)}}$:} Taking the
coefficients of the 
monomials $x^2y^{i-2}$  in~\eqref{rotational2} and the
monomial $x^3y^{i-1}$ in~\eqref{integral2}, we obtain the system :
\begin{align*}
&3b^{(3,i-2,0)}-(i-1)a^{(2,i-1,0)}=\left(\begin{array}{c} i \\2 \end{array}\right)(\lambda_1\cos(\frac{(i-2)\pi}{2})+\lambda_2\sin(\frac{(i-2)\pi}{2})),\\
&a^{(2,i-1,0)}+b^{(3,i-2,0)}=0,
\end{align*}
which has determinant $-i-2$. The system is completely determined and the solution for $a^{(2,i-1,0)}$ is 
\[
a^{(2,i-1,0)}=\displaystyle\frac{(-i^2+i)\left(\lambda_1\cos(\frac{(i-2)\pi}{2})+\lambda_2\sin(\frac{(i-2)\pi}{2}\right)}{2(i+2)}.
\]
\vspace*{0.2cm}

\noindent{\it -Computation of $\mathbf{b^{(1,i,0)}}$:} Taking the
coefficients of the monomial $y^{i}$  in~\eqref{rotational} and the monomial $xy^{i+1}$ in~\eqref{integral} yields: 
\begin{align*}
&b^{(1,i,0)} -(i+1)a^{(0,i+1,0)}=\lambda_1\cos(\frac{i\pi}{2})+\lambda_2\sin(\frac{i\pi}{2}),\\
&a^{(0,i+1,0)}+b^{(1,i,0)}=0,
\end{align*}
which has determinant $-i-2$. The system is completely determined and the solution for $b^{(1,i,0)}$ is given by
\[
b^{(1,i,0)}=\displaystyle\frac{\lambda_1 \cos \left(\frac{i\pi }{2}\right)+\lambda_2 \sin \left(\frac{i\pi }{2}\right)}{i+2}.
\]

\vspace*{0.2cm}
\noindent{\it -Computation of $\mathbf{c^{(1,i-1,1)}}$:} Taking the
coefficients of the  monomials $y^{i-1}z$ and $xy^{i-2}z$ in~\eqref{rotational2} and the monomial $xy^{i-1}z^2$ in~\eqref{integral2} we obtain the system:
\begin{align*}
&2a^{(0,i-1,2)}-c^{(1,i-1,1)}=i\left(\lambda_1\cos(\frac{i\pi}{2})+\lambda_2\sin(\frac{i\pi}{2})\right),\\
&-2b^{(1,i-2,2)}+(i-1)c^{(1,i-1,1)}=\left(\begin{array}{c}i\\2\end{array}\right)2\left(\lambda_1\cos(\frac{(i-2)\pi}{2})+\lambda_2\sin(\frac{(i-2)\pi}{2})\right),\\
&a^{(0,i-1,2)}+b^{(1,i-2,2)}-ic^{(1,i-1,1)}=0,
\end{align*}
which has determinant $2i$ and whose solution for $c^{(1,i-1,1)}$ is
\[
c^{(1,i-1,1)}=i\left(\lambda_1 \cos \left(\frac{i\pi }{2}\right)+\lambda_2 \sin \left(\frac{i\pi }{2}\right)\right).
\]
Substituting the variables by their values into Equation~\eqref{divergences2} yields
\begin{align*}
&\left(\displaystyle\frac{-i}{i+2}   -1\right)          \left(\lambda_1 \cos \left(\frac{(i-1)\pi}{2}\right)+ \lambda_2 \sin
   \left(\frac{(i-1)\pi}{2}\right)\right)=0,
\\
&\left(\displaystyle\frac{i^2 - i}{i + 2} + \frac{i}{2 + i} + i\right)\left(\lambda_1\cos(\frac{i\pi}{2})+\lambda_2\sin(\frac{i\pi}{2})\right)=0.
\end{align*}
Finally, arguing as in the proof of Proposition~\ref{con3} we obtain the desired result.

\qed
\section{A Strong Unique Continuation Principle}\label{app:strong_continuation}
In this appendix we recall a unique continuation theorem which implies that any $C^\infty$ Beltrami field with a zero of infinite order (i.e., the whole Taylor expansion at the zero point vanishes) is identically zero. This property is extensively used in the proof of our main theorems.

\begin{theorem}[see \cite{AroKrzySza}] Let $X$ be a $C^\infty$ vector field defined in a bounded domain $K\subset\mathbb R^3$ and satisfying 
\[
\|curl(X)\|_{L^\infty}+\|div(X)\|_{L^\infty} \leq c\|X\|_{L^\infty} 
\]
for some constant $c>0$. If $p\in K$ is a zero point of $X$ of infinite order, that is, all the derivatives $D^{\alpha}X(p)$ vanish for any multi-index
$\alpha$, then $X\equiv 0$ in $K$.
\end{theorem}
Notice that, in particular, the bound for the $L^\infty$ norm in this theorem holds for any Beltrami field, with $c$ the supremum of $f$ in the set $K$.

\section{An Open Problem}\label{app:numerical}

In this final appendix, we explore the question of whether the algebraic obstructions found in Theorem~\ref{th4} can be
applied to include non-vanishing terms $f_3$ and $f_5$ in the series
expansion of $f$ at $p$. The simplest case is when
$f_0\neq 0$, where the only term that we are not able to include in
our results (see Remark~\ref{remark}) is $f_3$. The algebraic obstructions found in this article have an irregular behavior when
applied to functions with non-trivial term $f_3$. 

Indeed, for functions $f$ where the spectrum of the Hessian is of the form
$\{\alpha,-\alpha,\beta\}$, the obstruction used in Theorem~\ref{th4} is not enough to obtain
the desired result: non-existence of non-trivial solutions. That is, the set of equations
\begin{align*}
  &curl(X_i)=0,\\
  &div(X_{i})=0,\\
  &\langle \nabla(\sigma_1 x^2+\sigma_2 y^2+\sigma_3 z^2), X_{i} \rangle=0,\\
&curl(X_{i+1})=f_0X_i,
\\
&div(X_{i+1})=0,
\\
&\langle \nabla(\sigma_1 x^2+\sigma_2 y^2+\sigma_3 z^2), X_{i+1} \rangle+\langle \nabla f_3,X_i\rangle=0,
\end{align*}
does not guarantee that $X_i\equiv 0$. For instance, consider
$$f=1+\frac12(x^2+y^2-z^2)+2xyz.$$ In that case $p=0$ is a non-degenerate
critical point where the function does not vanish. The eigenvalues of
the Hessian are $\{1,1,-1\}$, so the first possible non-vanishing term
is $X_1$. It can be seen that $X_1=\lambda_1X_1^1+\lambda_2X_1^2$,
where
\[
\begin{array}{lr}
X_1^1=\left(\begin{array}{c}
z \\0 \\x
\end{array}\right),
&
X_1^2=\left(\begin{array}{c}
z \\0 \\y
\end{array}\right).
\end{array}
\] 
If we impose the following equations on the term $X_2$ 
\begin{align*}
&curl(X_{2})=f_0X_1,
\\
&div(X_{2})=0,
\\
&\langle \nabla(x^2+ y^2-z^2), X_{2}
  \rangle+\langle  \nabla f_3, X_{i}\rangle=0,
\end{align*}
we can find the non trivial solutions
\[
\begin{array}{lr}
X_2=\left(\begin{array}{c}
x y\\ 0\\ -y z
\end{array}\right),
&
X_1=\left(\begin{array}{c}
-z\\0\\ -x
\end{array}\right)
\end{array}
\] 
which show that this obstruction is not enough to conclude $X_1 \equiv 0$ and therefore $X \equiv 0$. How many equations should one include to obtain the desired result should be the object of future research. This should be pursued alongside the application of more sophisticated algebraic techniques.

\end{appendices}

\nocite{*}
\bibliography{template}
\bibliographystyle{acm}

\end{document}